\documentclass[reqno,12pt]{article}
\usepackage{amsfonts,amssymb,latexsym,amsmath,amsthm,a4}

\newtheorem{theorem}{Theorem}[section]
\newtheorem{lemma}[theorem]{Lemma}

\theoremstyle{definition}
\newtheorem{remark}[theorem]{Remark}
\newtheorem{definition}[theorem]{Definition}
\newtheorem{example}[theorem]{Example}

\numberwithin{equation}{section}

\def\extp{$1.12_{\lambda,p}$}
\def\exst{$1.12_{\lambda^*,p}$}
\def\ve{\varepsilon}
\def\Ren{\mathbb{R}^n}
\def\Om{\Omega}
\def\b1z{B_1\setminus\{0\}}
\def\B1zc{\overline{B}_1\setminus\{0\}}
\def\W1p{{W^{1,p}(B_1)}}
\def\R{\mathbb{R}}
\def\aver{\operatornamewithlimits{{\int}\kern-11 true pt \diagup}}

\begin{document}

\title{Regularity of radial minimizers of reaction
equations involving the $p$-Laplacian}

\author{
\sc Xavier Cabr\'e\\
{\sc\footnotesize ICREA and Universitat Polit\`ecnica de Catalunya}\\
{\sc \footnotesize Departament de Matem\`{a}tica  Aplicada I}\\
{\sc \footnotesize Diagonal 647, 08028 Barcelona, Spain}\\
{\tt\footnotesize xavier.cabre@upc.edu}
\and
\sc Antonio Capella\\
{\sc\footnotesize Universit\"at Bonn}\\
{\sc \footnotesize Institut f\"ur Angewandte Mathematik}\\
{\sc \footnotesize Wegelerstr. 10, 53115 Bonn, Germany}\\
{\tt\footnotesize capella@iam.uni-bonn.de}
\and
\sc Manel Sanch\'on\\
{\sc\footnotesize Universitat de Barcelona}\\
{\sc \footnotesize Departament de Matem\`atica Aplicada i An\`alisi}\\
{\sc \footnotesize Gran Via 585, 08007 Barcelona, Spain}\\
{\tt\footnotesize msanchon@maia.ub.es} }
  

\maketitle

\begin{abstract}
We consider semi-stable, radially symmetric, and decreasing solutions
of $-\Delta_p u=g(u)$ in the unit ball of $\R^n$, where $p>1$,
$\Delta_p$ is the $p$-Laplace operator, and $g$ is a locally
Lipschitz function. For this class of radial solutions, which
includes local minimizers, we establish pointwise, $L^q$, and
$W^{1,q}$ estimates which are optimal and do not depend on the
specific nonlinearity $g$. Among other results, we prove that every
radially decreasing and semi-stable solution $u$ belonging to
$W^{1,p}(B_1)$ is bounded whenever $n<p+4p/(p-1)$.

Under standard assumptions on the nonlinearity $g(u)=\lambda f(u)$,
where $\lambda>0$ is a parameter, it is proved that the
corresponding extremal solution $u^*$ is semi-stable, and hence, it
enjoys the regularity stated in our main result.
\end{abstract}

\section{Introduction}
\label{section1}

This article is concerned with reaction-diffusion equations
$-\Delta_p u=g(u)$
involving the $p$-Laplace operator. We consider a class of
radially symmetric solutions: those solutions which are decreasing
and semi-stable. This type of solutions include local minimizers,
minimal solutions, extremal solutions, and also certain solutions
found between a sub and a supersolution. We establish sharp
pointwise, $L^q$, and $W^{1,q}$ estimates for solutions in this
class. In particular we show that semi-stable radial solutions enjoy
better regularity properties than general radial solutions. 
In addition, our 
results do not depend on the specific form of the nonlinearity in
the reaction term. More precisely, our pointwise and $L^q$
estimates hold for every locally Lipschitz nonlinearity $g$, while
the $W^{1,q}$ estimates hold for nonnegative $g$. All the results
obtained in this paper were obtained by the first two authors in
\cite{CC05} for the Laplace operator, namely $p=2$.

We consider radially symmetric and decreasing solutions $u\in\W1p$ of
\begin{equation}\label{problem1}
-\textrm{div}(|\nabla u|^{p-2}\nabla u) = g(u)\quad \textrm{in
}B_{1}\setminus\{0\},
\end{equation}
where $p>1$, $B_1$ is the unit ball of $\Ren$, and
$g:\R\longrightarrow \R$ is a locally Lipschitz function. By a
radially decreasing function $u$ we mean a function $u$ such that
$u=u(r)$ and $u_r(r)=(du/dr)(r)<0$ for all $r\in(0,1)$, where
$r=|x|$ and $u_r$ denotes the radial derivative.

We do not assume $u(1)=0$ or any other boundary condition. Nevertheless, 
since $u\in \W1p$ is radial, we have $u\in L^\infty_{\rm loc}
(\B1zc)$ by the Sobolev embedding in one dimension. 
Hence, using known regularity results for
degenerate elliptic equations (see \cite{Lie88}), we have that
in fact $u\in C^{1,\beta}_{\rm loc}(\overline{B}_1\setminus\{0\})$
for some $\beta\in(0,1)$. However, $u$ may be unbounded at the
origin. Our main result establishes estimates for semi-stable
solutions in the whole ball $B_1$. Thus, this may be regarded as a result
on removable singularities.

Consider the energy functional
\begin{equation} \label{functional}
E_\Om(u):=\frac{1}{p}\int_\Om |\nabla u|^p\, dx-\int_\Om G(u)\, dx,
\end{equation}
where $G'=g$ and $\Om$ is a smooth bounded domain of $\Ren$. We say
that a radially decreasing function $u\in \W1p$ is a \emph{radial
local minimizer} of \eqref{functional} with $\Omega=B_1$ if for
every $\delta>0$ there exists $\varepsilon_\delta>0$ such that
$$
E_{B_1\setminus\overline{B}_\delta}(u)\leq
E_{B_1\setminus\overline{B}_\delta}(u+\xi)
$$
for all radial functions $\xi\in C^1_c(B_1\setminus\overline{B}_\delta)$ 
(i.e., with compact support in $B_1\setminus \overline{B}_\delta$) 
satisfying $\|\xi\|_{C^1}\leq \varepsilon_\delta$. 
Note that every radial local minimizer $u$ is a
solution of \eqref{problem1}. Moreover, $u$ is semi-stable in the
sense that the second variation of the energy at $u$ is nonnegative.
The following definition makes this precise ---see Remark~\ref{rem1} 
for more comments on this, as well as \cite{CS05} for the corresponding
setting in the nonradial case. 

\begin{definition}\label{def:stab}
Let $u\in \W1p$ be a radially symmetric solution in $B_1\setminus\{0\}$
of \eqref{problem1}
such that $u_r(r)<0$ for all $r\in (0,1)$. We say that $u$ is
\emph{semi-stable} if and only if the second variation of energy $Q$ 
at $u$ satisfies
\begin{equation}\label{semistab}
Q(\xi):=\int_{B_1}\Big\{(p-1)|u_r|^{p-2}|\xi_r|^2-g'(u)\xi^2\Big\}\,
dx\geq 0,
\end{equation}
for every radially symmetric function $\xi\in C^1_c(\b1z)$.
\end{definition}
As we will see later, the class of semi-stable solutions includes not
only local minimizers but also the minimal and extremal solutions of 
the classical problem (\extp) below, which motivated our work.

To state our estimates, we define exponents $q_k$, for $k=0,1$, by
\begin{equation}\label{q0_and_q1}
\left\{ \begin{array}{ll} \displaystyle
\frac{1}{q_{_k}}:=\frac{1}{p}-\frac{2}{np}
\sqrt{\frac{n-1}{p-1}}+\frac{k-1}{n} -\frac{2}{np}
&\text{ for }\displaystyle n \ge p+\frac{4p}{p-1}\\
q_{_k} :=+\infty &\text{ for }\displaystyle n < p+\frac{4p}{p-1}.
\end{array}\right.
\end{equation}
One can easily check that $p<q_{_k}\leq +\infty$ in all cases. In
addition,
\begin{equation}\label{rmk:cot_np}
\frac{n-p}{2} \ge 1+\sqrt{\frac{n-1}{p-1}} \qquad
\text{ if and only if } \qquad n\ge p+\frac{4p}{p-1}.
\end{equation}
Hence, both exponents $q_0$ and
$q_1$ are well defined. Moreover, $q_0$ and $q_1$ are finite if
$n>p+4p/(p-1)$.

Our main theorem states sharp regularity results and pointwise
estimates for every semi-stable radially decreasing solution $u$ of
\eqref{problem1}. It also establishes estimates for its radial
derivative $u_r$. Our result extends Theorem~1.5 and part of
Theorem~1.8 of \cite{CC05} (in which $p=2$) to the general case
$p>1$. In particular, if we set $p=2$ in the following theorem, one
recovers the results of \cite{CC05} for the Laplace operator.
The statements in the following theorem are optimal in several 
respects discussed below, and are still open in general domains
even for $p=2$. Indeed, in the nonradial case, they are known to 
hold only for $g(u)=\lambda e^u$, $g(u)=\lambda(1+u)^m$ with 
$m>p-1$, and nonlinearities $g$ close to them in certain senses
---see the comments following Theorem~\ref{thm2} for more details. 

\begin{theorem}\label{thm1}
Let $g$ be a locally Lipschitz function and $u\in W^{1,p}(B_1)$ be a
semi-stable radial solution in $B_1\setminus\{0\}$ of \eqref{problem1} satisfying
$u_r(r)<0$ for all $r\in (0,1)$. Then:

{\rm (a)} If $n<p+4p/(p-1)$ then $u\in L^\infty(B_1)$. Moreover,
%
$$
 \Vert u\Vert_{L^\infty(B_1)} \leq C_{n,p}\Vert u
\Vert_{W^{1,p}(B_1)},
$$
%
where $C_{n,p}$ is a constant depending only on $n$ and $p$.

{\rm (b)} If $n=p+4p/(p-1)$ then $u\in L^q(B_1)$ for all
$q<+\infty$. Moreover,
\begin{equation}\label{eq:logpbd}
|u(r)| \leq C_{p} \Vert u \Vert_{W^{1,p}(B_1)} (|\log r
|+1)\quad\text{ in }B_1,
\end{equation}
where $C_{p}$ is a constant depending only on $p$.

{\rm (c)} If $n> p+ 4p/(p-1)$ and $q<q_0$, then $u\in L^q(B_1)$ and
%
$$
\Vert u \Vert_{L^q(B_1)}\leq C_{n,p,q} \Vert u \Vert_{W^{1,p}(B_1)},
$$
%
where $C_{n,p,q}$ is a constant depending only on $n$, $p$, and $q$.
Moreover,
\begin{equation}\label{eq:estLqpws}
|u(r)|\leq C_{n,p}\Vert u \Vert_{W^{1,p}(B_1)}
r^{-\frac{1}{p}\left(n-2\sqrt{\textstyle\frac{n-1}{p-1}}-p-2\right)} (|\log
r|^\frac{1}{p}+1)\text{ in }B_1,
\end{equation}
where $C_{n,p}$ is a constant depending only on $n$ and $p$.

{\rm (d)} Assume that $g$ is nonnegative. Then:

{\rm (d1)} We have 
\begin{equation}\label{eq:estGLq}
\Vert \nabla u\Vert_{L^p(B_1)}
\le C_{n,p}\left\{\Vert (u-u(1))^{p-1} \Vert_{L^1(B_1)}^\frac{1}{p-1}
+\Vert g(u)\Vert_{L^1(B_1)}^{\frac{1}{p-1}} \right\}
\end{equation}
for some constant $C_{n,p}$ depending only on $n$ and $p$.

{\rm (d2)} $u\in W^{1,q}(B_1)$ for every $q<q_1$, and 
\begin{equation}\label{eq:Wbd}
\Vert u\Vert_{W^{1,q}(B_{1})}\le C\, \qquad\text{ if } q<q_1,
\end{equation}
where $C$ is a constant depending only on $n$, $p$, $q$, and on
upper bounds for $\Vert u\Vert_{L^1(B_1)}$ and  $g$.  

{\rm (d3)} If $n\ge p+4p/(p-1)$ then
\begin{equation}\label{eq:pdWdb}
|u_r(r)|\leq C_{n,p} \Vert u \Vert_{W^{1,p}(B_1)}
r^{-\frac{1}{p}\left(n-2\sqrt{\textstyle\frac{n-1}{p-1}}-2\right)}|\log
r|^{\frac{1}{p}}\text{ in }B_{1/4},
\end{equation}
where $C_{n,p}$ is a constant depending only on $n$ and $p$.
\end{theorem}

Examples \ref{exp:plap1} and \ref{exp:plap2} below on the
exponential and power nonlinearities show the sharpness 
of the previous regularity results.
Indeed, the functions 
\begin{equation}\label{eq:1.15-1}
u(r)=-p\log r\quad\text{and}\quad u(r)=r^{-\frac{p}{m-(p-1)}}-1
\end{equation}
are solutions of \eqref{problem1} for, respectively,
$g(u)=\lambda^*e^u$ and $g(u)=\lambda^*(1+u)^m$ with certain values
of $\lambda^*=\lambda^*(n,p,m)$. Using weighted Hardy inequalities, 
one can easily find the ranges of values of $n$, $p$, and $m$ for 
which these solutions are $\W1p$ semi-stable solutions (see \cite{CS05}
and Examples \ref{exp:plap1} and \ref{exp:plap2} below for more details).
The ranges obtained in this way for the exponential nonlinearity show 
the sharpness of condition $n<p+4p/(p-1)$ in Theorem~\ref{thm1}(a), 
as well as the optimality of the pointwise estimate \eqref{eq:logpbd}
in part (b). The ranges obtained for the power nonlinearities give  
the optimality of the exponents $q_0$ and $q_1$ in parts (c) and (d), 
as well as the sharpness of the pointwise bounds \eqref{eq:estLqpws}
and \eqref{eq:pdWdb} except for the factors $|\log r|^{1/p}$ in them.


The proof of Theorem \ref{thm1} was inspired in the proof of Simons
theorem on the nonexistence of singular minimal cones in
$\mathbb{R}^n$ for $n\le 7$ (see \cite{CC06} for further details).
The main point of the proof is to obtain the following key estimate
$$
\int_{B_1}|u_r|^pr^{-2\alpha}\, dx \leq
\frac{C_{n,p}}{(n-1)-(\alpha-1)^2(p-1)} \Vert
\nabla u\Vert_{L^{p}(B_1)}^p
$$
---from which most of our results follow--- for those exponents
$\alpha$ making the denominator in the above right-hand side positive
(see Lemma~\ref{lem2:plap}).
This estimate is established by taking $\xi=u_r\eta$ as a
test function on the semi-stability condition \eqref{semistab},
being $\eta$ essentially a negative power of $r$.
Lemma \ref{lem2:1} shows that, with this choice of $\xi$ in 
\eqref{semistab}, the term $g'(u)$
in  \eqref{semistab} disappears. 
This is the reason why our main estimates do not depend on $g$.


As an application of Theorem~\ref{thm1} and the ideas behind its
proof, we consider the following problem
\stepcounter{equation}
$$
\left\{
\begin{array}{rcll}
-\text{div}(|\nabla u|^{p-2}\nabla u)&=&\lambda f(u)&
\textrm{ in }B_1,\\
u& > & 0 & \textrm{ in }B_1,\\
u& = & 0 & \textrm{ on
}\partial B_1,
\end{array}
\right.\eqno{(1.12_{\lambda,p})}
$$
where $\lambda>0$, and $f$ is an increasing $C^1$ function with
$f(0)>0$ and
\begin{equation}\label{superlinear}
\lim_{t\rightarrow+\infty}\frac{f(t)}{t^{p-1}}=+\infty.
\end{equation}

Problem (\extp) is studied in \cite{CS05} for general
smooth bounded domains $\Omega$ of $\R^n$. It is proved there that
there exists a parameter $\lambda^*\in (0,\infty)$ such that if $0<\lambda
<\lambda^*$ then (\extp) admits a minimal solution
$u_\lambda\in C^1 (\overline{\Omega})$. Here minimal means smaller
than any other supersolution of the problem. Moreover, for
$\lambda>\lambda^*$ problem (\extp) admits no regular
solution. We also have that for every $0< \lambda < \lambda^*$ the
minimal solution $u_\lambda$ is semi-stable (in the sense of
Definition~\ref{def:stab} when $\Omega=B_1$).
Consider the increasing limit
\begin{equation}\label{extremal}
u^*:=\lim_{\lambda\uparrow \lambda^*}u_\lambda.
\end{equation}
In contrast with the case $p=2$ involving the Laplacian, 
it is not always clear that the
limit $u^*$ is a weak solution of (\extp) for
$\lambda=\lambda^*$. When one can establish that $u^*$ is 
a weak solution (this may depend on the assumptions on $n$, 
$p$, $\Omega$, and $f$), it is called the \emph{extremal solution}.

In our next result we prove that, for $\Omega=B_1$, the limit 
$u^*$ is actually a semi-stable radially decreasing energy solution of
(\exst). As a consequence, it enjoys the same
regularity properties as the ones stated in Theorem~\ref{thm1}. In
particular, we obtain that the extremal solution $u^*$ is bounded if
$ n< p+4p/(p-1)$. As we will see below, the extremal solution is 
unbounded if $n\geq p+4p/(p-1)$ when $f(u)=e^u$. Therefore, this 
range of dimensions is optimal. The optimal dimension ensuring
the boundedness of the extremal solution  in general domains
remains unknown for general 
nonlinearities $f$ under the assumptions above (and even 
with the additional hypothesis that $f$ is convex).

\begin{theorem}\label{thm2}
Let $f$ be a positive and increasing $C^1$ function in $[0,+\infty)$ 
satisfying \eqref{superlinear}. For $\lambda\in (0,\lambda^*)$  
let $u_\lambda$ be the minimal solution of {\rm (\extp)},
and let $u^*$ be defined by \eqref{extremal}.

Then, 
\begin{equation}\label{eq:unifbd}
\Vert u_\lambda\Vert_\W1p+\Vert f(u_\lambda)\Vert_{L^1(B_1)}\le C
\end{equation}
for all $\lambda\in (0,\lambda^*)$ and some constant $C$ independent
of $\lambda$. Moreover, $u^*\in\W1p$ and $u^*$ is a semi-stable 
radially decreasing energy solution of {\rm (\exst)}. 
As a consequence, $u^*$ has the regularity stated in 
Theorem~{\rm \ref{thm1}}. In particular, 
\begin{equation}\label{ext:boundedness}
u^*\in L^\infty(B_1)\quad\text{if}\quad n<p+\frac{4p}{p-1}.
\end{equation}
\end{theorem}

As mentioned in the statement of the theorem, for problem 
{\rm (\extp)} we consider \emph{energy solutions}. 
That is, nonnegative functions  $u\in W^{1,p}_0(B_1)$ such 
that $f(u)\in L^1(B_1)$ and
\begin{equation}\label{distributional}
\int_{B_1}|\nabla u|^{p-2}\nabla u\cdot\nabla\varphi\,
dx=\int_{B_1}\lambda f(u)\varphi\, dx\quad \textrm{for all
}\varphi\in C_c^1(B_1).
\end{equation}
Note that, by a standard density argument, \eqref{distributional}
holds for every $\varphi\in W^{1,p}_0(B_1)$.

The literature studying minimal and extremal solutions in general
domains of $\R^n$ is extensive. Crandall and Rabinowitz \cite{CranRa75} 
and Mignot and Puel \cite{MP80} considered the case of the exponential 
and power nonlinearities in general domains for $p=2$. They proved that $u^*$
is an energy solution for every dimension and that it is bounded in
some range of dimensions. These results in general domains were extended 
for every $p>1$ by Garc\'{\i}a-Azorero, Peral, and 
Puel~\cite{GarPe92,GarPePu94} for the exponential case, and  by 
Ferrero~\cite{Fe04} and Cabr\'e and Sanch\'on~\cite{CS05} for power type 
nonlinearities.

An optimal regularity result for the extremal solution is not known for
general nonlinearities and general domains, even for $p=2$,
with the exception of the radial case $\Omega=B_1$.
In \cite{CC05}, Cabr\'e and Capella obtained the optimal regularity
of the extremal solution for every locally Lipschitz nonlinearity
$f$ when the domain is a ball and $p=2$. In particular, they proved
that the extremal solution is bounded if $n\leq 9$. In this paper, we 
obtain an analogue optimal radial result for all $p>1$.

In the nonradial case, \cite{ND} and \cite{2a} contain the best results
for $p=2$. In \cite{ND}, Nedev proves for $p=2$ that the extremal solution
is bounded if $n\le 3$ ---i.e., $n<4$ in \eqref{boundedness:Nedev}---  
 whenever $f$ is an increasing and convex function such that $f(0)>0$ and 
\eqref{superlinear} holds. The same result up to dimension $n\le 4$ is proved
by the first author in \cite{2a} for $p=2$ and $\Omega$ convex, without
the convexity assumption on $f$. On the other hand, the third author extends 
in \cite{S} the work of Nedev \cite{ND} to the case $p>2$ and establishes
\begin{equation}\label{boundedness:Nedev}
u^*\in L^\infty (\Omega) \quad\text{if}\quad n<p+\frac{p}{p-1}.
\end{equation}
Note that for $p=2$, we recover Nedev's condition $n<4$. Note also the
existing gap between the dimension in \eqref{boundedness:Nedev} for 
general $\Omega$ and the ones in our result for $\Omega=B_1$ 
---note the presence of the factor $4$ in \eqref{ext:boundedness}
in contrast with \eqref{boundedness:Nedev}. 

We refer to \cite{8bis,2b} for surveys on minimal and extremal 
solutions and to \cite{BV,46,47,48,EE05,ESP,GG,JL,S2} for other
interesting results in the topic of extremal solutions ---also for
problems involving other operators, such as the bilaplacian, or
fractional Laplacians related to boundary reactions.



The key point in the proof of Theorem~\ref{thm2} is to establish 
the uniform estimate \eqref{eq:unifbd} on the $W^{1,p}$-norm and the
$L^1(B_1)$-norm of $u_\lambda$ and $f(u_\lambda)$, respectively.
This is accomplished using the ``superlinearity'' hypothesis 
\eqref{superlinear} on $f$, together with the radially decreasing 
character of the solutions. The $W^{1,p}$ bound is then used to 
show that the limit $u^*$ is an energy solution of 
(\exst).
Using again that the $u_\lambda$ are semi-stable, we prove that $u^*$ is
semi-stable.


In \cite{CS05}, the semi-stability property of the extremal solution
$u^*$ was proved in the case $p\geq 2$ for general domains in $\R^n$
and nonlinearities having its growth comparable to a power. This
property was obtained taking the limit as $\lambda\to\lambda^*$
in the semi-stability condition for the minimal solutions $u_\lambda$.
In the case $1<p<2$ this argument does not apply for general
domains, since we have no control on the set where the gradient
$\nabla u^*$ vanishes. However, in the radial case
this set is the origin, and thus we obtain here the result of \cite{CS05}
also for $1<p<2$.


Next, we introduce two explicit examples which show the
optimality of Theorems \ref{thm1} and \ref{thm2}.

\begin{example}\label{exp:plap1}In \cite{GarPe92,GarPePu94}
Garc\'{\i}a-Azorero, Peral, and Puel
considered problem (\extp) for $f(u)=e^u$ and a general
bounded domain $\Omega$. They proved
that if $n<p+4p/(p-1)$ then the extremal solution $u^*$ is bounded.
On the other hand, if $n\geq p+4p/(p-1)$ and $\Omega=B_1$ they show 
that
\begin{equation}\label{eq:1.14-1}
u^*(r)=-p\log r\quad\text{and}\quad \lambda^*=p^{p-1}(n-p),
\end{equation}
and hence, that the extremal solution is
unbounded. This example shows that the conditions on the dimension
$n$ in terms of $p$ in Theorem~\ref{thm1}(a),(b) and the pointwise
estimate \eqref{eq:logpbd} are optimal.

Indeed, as mentioned before, it is easy to check that for the previous
range of dimensions, $u^*$ as in \eqref{eq:1.14-1} is a singular $W^{1,p}$ 
semi-stable solution ---independently of the more precise fact 
of actually being the extremal solution.  This remark also applies to 
the following example.
\end{example}

\begin{example}\label{exp:plap2}
Consider the power nonlinearity $f(u)=(1+u)^m$, with $m>p-1$,
which is studied in \cite{CS05} and \cite{Fe04} for general domains. 
Let us define the following critical exponent 
$$
m_{cs}(p):= \left\{
\begin{array}{ll}
\displaystyle \frac{(p-1)n-2\sqrt{(p-1)(n-1)}+2-p}
{n-(p+2)-2\sqrt{\frac{n-1}{p-1}}} &\textrm{if }\displaystyle
n>p+\frac{4p}{p-1},\\
+\infty&\textrm{if }\displaystyle n\leq p+\frac{4p}{p-1}.\\
\end{array}
\right.
$$
The results of \cite{CS05,Fe04} show that if $m<m_{cs}(p)$ then 
the extremal solution $u^*$ of (\extp) is bounded. 
On the other hand, if $\Omega=B_1$ and  $m\geq m_{cs}(p)$ then
$$
u^*(r)=r^\frac{-p}{m-(p-1)}-1\quad \textrm{and}\quad
\lambda^*=\left(\frac{p}{m-(p-1)}\right)^{p-1}
\left(n-\frac{mp}{m-(p-1)}\right).
$$
We see that $u^*\in L^q(B_1)$ if and only if $1\leq q<
n(m-(p-1))/p$. Now, if we let $m=m_{cs}(n,p)$ then we find that
$n(m_{cs}(n,p)-(p-1))/p=q_0$ and this proves the sharpness of
Theorem~\ref{thm1}(c).

Observe that in the same case $m=m_{cs}(n,p)$ we have
$$
u^*(r)=r^{-\frac{1}{p}\left(n-2\sqrt{\frac{n-1}{p-1}}-p-2\right)}-1,
$$
and, as in the case $p=2$, this differs from the pointwise bound
\eqref{eq:estLqpws} for the factor $|\log r|^{1/p}$. Recently,
Villegas \cite{V07b} has proved that, for $p=2$, the factor $|\log
r|^{1/2}$ in \eqref{eq:estLqpws} can be removed. See also
\cite{V07a} for a related improvement of a radial Liouville 
theorem of \cite{CC}.
\end{example}

\begin{remark}
The critical dimension for the boundedness of semi-stable solutions
determined by $n=p+4p/(p-1)$ tends to $+\infty$ as $p\to 1$ and as
$p\to+\infty$. The smallest critical dimension corresponds to $p=3$,
for which $n=9$. In the cases $p=2$ and $p=5$ we find the value
$n=10$. Observe that case (b) in Theorem~\ref{thm2}, $n=p+4p/(p-1)$,
applies only for those $p>1$ such that $p+4p/(p-1)$ is an integer.
Note also that the critical dimension for the boundedness of all
semi-stable solutions is bigger than the corresponding dimension for 
general solutions, which is $n<p$.
Therefore, semi-stable solutions enjoy more regularity than general
solutions.
\end{remark}

\begin{remark}\label{rem1}
The following two comments concern the notion of semi-stability of 
solutions. See \cite{CS05} for a more general nonradial setting.

(i) If $u$ is a radially decreasing local minimizer as defined after 
\eqref{functional}, we claim that $u$ is a
semi-stable solution. Therefore, Theorem~\ref{thm1} applies to every
radially decreasing local minimizer $u\in W^{1,p}(B_1)$. To check
this, note that the second variation of the energy functional
\eqref{functional} in $B_1$ at a radially decreasing solution $u$ is
given by
\begin{equation}\label{second_var}
\int_{B_1}\left\{|\nabla u|^{p-2} \left((p-2)\left[\frac{\nabla
u}{|\nabla u|}
\cdot\nabla\xi\right]^2+|\nabla\xi|^2\right)-g'(u)\xi^2\right\}\, dx,
\end{equation}
for every perturbation $\xi\in C^1$ with compact support in
$B_1\setminus \{0\}$ (not necessarily radially symmetric). It
follows that, if $u$ is a radial local minimizer, then
\eqref{second_var} is nonnegative for every radial $\xi\in C^1_c
(B_1\setminus\{0\})$. Therefore, $u$ is semi-stable in the sense of
Definition~\ref{def:stab} since for radial perturbations $\xi$,
\eqref{second_var} reduces to \eqref{semistab}.

(ii) If $u$ is a radially decreasing semi-stable solution of
\eqref{problem1} in the sense of Definition~\ref{def:stab}, we claim 
that \eqref{second_var} is nonnegative for every nonradial $\xi\in C^1$
with compact support in $B_1\setminus\{0\}$. That is, the second variation
of energy at $u$ is nonnegative not only for radial perturbations but 
also for nonradial ones. Indeed, for $x\in \mathbb{R}^n$ we
set $x=r\theta$, where $r=|x|\geq 0$, and $\theta=x/r\in
\partial B_1=\{\theta\in\mathbb{R}^n\, :\, |\theta|=1\}$.
Let $\xi\in C^1_c(B_1\setminus\{0\})$ (not necessarily radial). We
consider the spherical averages of $\xi^2$ and define the following
radial function:
$$
\varphi^2(r):=\frac{1}{|\partial B_1|} \int_{\partial
B_1}\xi^2(r\theta)\ d\theta=\aver_{\partial B_1}\xi^2(r\theta)\
d\theta.
$$
Differentiating the last expression with respect to $r$ and using
the Cauchy-Schwarz inequality, we find
$$
\varphi_r^2(r) \leq \aver_{\partial B_1}\left[\frac{x}{r}
\cdot\nabla\xi(r\theta)\right]^2d\theta =\aver_{\partial
B_1}|\xi_r(r\theta)|^2d\theta.
$$
Finally, after some straightforward computations, and using 
the semi-stability of $u$ in the sense of Definition~\ref{def:stab} 
for the radial test function $\varphi=\varphi(r)$ and that 
$\xi_r^2\leq |\nabla\xi|^2$ in  \eqref{second_var}, the claim follows.
\end{remark}


The paper is organized as follows. In Section \ref{section2} we
prove the main estimates on semi-stable stable solutions needed in the rest
of the paper. In Section \ref{section3} we prove Theorem~\ref{thm1}, while
Section~\ref{section4} is devoted to prove Theorem~\ref{thm2}.

\section{Estimates for semi-stable solutions}\label{section2}
This section is devoted to show  Lemma~\ref{lem2:plap}, which contains the key
estimate used in the proof of  Theorem~\ref{thm1}.

We begin with the following remark on radial functions and
radially decreasing solutions of \eqref{problem1}.
\begin{remark}\label{rem2:1}
Given $\delta\in (0,1)$, every radial function in $\W1p$ also belongs
(as a function of $r=|x|$) to the Sobolev space $W^{1,p}(\delta,1)$ 
in one dimension. As a consequence, by the Sobolev embedding in one 
dimension, $u$ is a continuous function of $r\in [\delta,1]$, 
\begin{equation}     
\label{sobolev_rad}
|u(1)|\le C_{n,p}\Vert u\Vert_\W1p\quad\text{and}\quad 
\Vert u\Vert_{L^{\infty}(B_1\setminus\overline{B}_\delta)}
\le C_{n,p,\delta}\Vert u\Vert_{\W1p}
\end{equation}
for some constant $C_{n,p}$ (respectively, $C_{n,p,\delta}$) depending 
only on $n$ and $p$ (respectively, on  $n$, $p$, and $\delta$).

Now, let $u\in\W1p$ be a radial solution to \eqref{problem1} such 
that $u_r(r)< 0$ for all $r\in (0,1)$. By \eqref{sobolev_rad}, the 
right-hand side of \eqref{problem1} is bounded away from the origin.
Hence, by standard regularity theory for the $p$-Laplacian, 
$u\in C^{1,\beta}_\text{loc}(\overline{B}_1\setminus\{0\})$ for some
$\beta\in (0,1)$. In particular, $u_r$ is a continuous function in
$B_1\setminus\{0\}$ which does not vanish and thus \eqref{problem1}
is uniformly elliptic in compact sets of $B_1\setminus\{0\}$. 
It follows that $u\in C^{2,\beta}_\text{loc}(B_1\setminus\{0\})$.

Since $u_r<0$ in $B_1\setminus\{0\}$, 
we can write \eqref{problem1} in radial coordinates as
\begin{equation}\label{problem1:radial}
-r^{1-n}\partial_r\left(r^{n-1}|u_r|^{p-2}u_r\right)=g(u)
\quad\text{ for }r\in (0,1),
\end{equation}
and also in the form
\begin{equation}\label{problem1:radial2}
-(p-1)|u_r|^{p-2}\partial_ru_r-\frac{n-1}{r}|u_r|^{p-2}u_r=g(u) \quad\text{ for
} r\in (0,1).
\end{equation}
\end{remark}

In order to show our main estimate (stated in Lemma~\ref{lem2:plap}) 
we need a preliminary result
on the form of the second variation of the energy for radial
solutions. This lemma was inspired in the proof of Simons theorem on
the nonexistence of singular minimal cones in $\mathbb{R}^n$ for
$n\le 7$ (see \cite{CC,CC06} for further details).

\begin{lemma}\label{lem2:1}
Let $u\in W^{1,p}(B_1)$ be a radial solution in $B_1\setminus\{0\}$ 
of \eqref{problem1} satisfying
$u_r(r)<0$ for all $r\in (0,1)$ and $Q$ be the quadratic form defined in 
\eqref{semistab}. Then,
\begin{equation}\label{lalala0}
Q(u_r\eta)=\int_{B_1}|u_r|^p\left\{(p-1)|\eta_r|^2
-\frac{n-1}{r^2}\eta^2 \right\}\, dx,
\end{equation}
for every radial function $\eta\in C^1_c(B_1\setminus\{0\})$, that is,
with compact support in $B_1\setminus\{0\}$.
\end{lemma}

Note that expression \eqref{lalala0} for the quadratic form $Q$
does not contain any reference to the nonlinearity $g$. This is the
reason why our estimates do not depend on the nonlinearity $g$.

\begin{proof}[Proof of Lemma {\rm \ref{lem2:1}}]
By Remark~\ref{rem2:1}, $u\in C^{2,\beta}_{{\rm loc}}(B_1\setminus\{0\})$ for
some $\beta\in(0,1)$ and $u$ satisfies equation \eqref{problem1:radial2}.
Let us set $H(u)=(p-1)|u_r|^{p-2}$. Let $\eta\in C^1_c(B_1\setminus\{0\})$ 
be a radial function with compact support in $B_1\setminus\{0\}$ and 
$c\in C^1(B_1\setminus\{0\})$ be a radial function. Take $\xi=c\eta\in
C^1_c(B_1\setminus\{0\})$ in \eqref{semistab} to obtain
\begin{eqnarray}
Q(c\eta) & = &\int_{B_1}\left\{H(u)|\nabla(c\eta)|^2-g'(u)c^2\eta^2
\right\} dx
\nonumber\\
&=&\int_{B_1}\left\{H(u)(c^2|\nabla\eta|^2+\nabla\eta^2\cdot c\nabla
c +\eta^2|\nabla c|^2) -g'(u)c^2\eta^2 \right\} dx \qquad
\nonumber\\ & = &\int_{B_1}\left\{H(u)c^2|\nabla\eta|^2
+H(u)\nabla(\eta^2c)\cdot\nabla c-g'(u)c^2\eta^2 \right\} dx.
\label{eq:6-08}
\end{eqnarray}

Next, we multiply \eqref{problem1:radial2} by $\partial_r(\eta^2u_r
r^{n-1})$, we integrate in $r$ from $0$ to $1$ and use integration 
by parts to obtain
$$
0\ = \ \int_0^1\!\!\!\partial_r(\eta^2u_r
r^{n-1})H(u)\partial_ru_rdr
-\int_0^1\!\!\!\partial_r\left\{\frac{n-1}{r}|u_r|^{p-2}u_r
+g(u)\right\}\eta^2u_r
r^{n-1}dr.
$$
Using $\partial_r(|u_r|^{p-2}u_r)=H(u)u_{rr}$ and 
$r^{n-1}dr=|\partial B_1|^{-1}dx$, we deduce
$$
0\ = \ \int_{B_1}
H(u)\partial_r(\eta^2u_r)\partial_ru_r\, dx-
\int_{B_1}\left\{g'(u)\eta^2u_r^2
-\frac{n-1}{r^2}|u_r|^{p-2}u_r^2\eta^2\right\}\, dx.
$$
Therefore
\begin{equation}\label{eq:6-10}
\int_{B_1}\Big\{H(u)\nabla(\eta^2u_r)\cdot\nabla u_r
-g'(u)\eta^2u_r^2\Big\}\, dx=-\int_{B_1}\frac{n-1}{r^2}\eta^2|u_r|^p\,
dx.
\end{equation}
Taking $c=u_r$ in \eqref{eq:6-08} and using \eqref{eq:6-10}, we
obtain \eqref{lalala0}.
\end{proof}


Now, we use Lemma~\ref{lem2:1} and the semi-stability
assumption to establish the following result. It is an
$L^p$ estimate for $u_rr^{-2\alpha/p}$ in $B_{1}$, for certain
positive exponent $\alpha$ depending on $n$ and $p$, in terms of the
$W^{1,p}$ norm of $u$. As said before, this is the key estimate
in the proof of Theorem~\ref{thm1}. Here as in the rest of this 
section, we assume that $n\ge p$. Note that when $n<p$, we have $\W1p
\subset L^\infty(B_1)$ and hence solutions are bounded. 

\begin{lemma}\label{lem2:plap}
Assume $n\geq p$. Let $u\in W^{1,p}(B_1)$ be a semi-stable radial
solution  in $B_1\setminus\{0\}$  
of \eqref{problem1} satisfying $u_r(r)<0$ for $r\in (0,1)$.
Let $\alpha$ satisfy
\begin{equation}\label{eq:alpha}
1\leq\alpha < 1+\sqrt\frac{n-1}{p-1}.
\end{equation}
Then
\begin{equation}\label{eq:6-11}
\int_{B_{1}}|u_r|^pr^{-2\alpha}\, dx \leq
\frac{C_{n,p}}{(n-1)-(\alpha-1)^2(p-1)} 
\Vert\nabla u\Vert_{L^p(B_1)}^p,
\end{equation}
where $C_{n,p}$ is a constant depending only on $n$ and $p$.
\end{lemma}

\begin{proof}
By the semi-stability of $u$ and Lemma~\ref{lem2:1} applied
with $\eta$ replaced by $r\eta$, we have that
\begin{equation}\label{stab:cond}
(n-1)\int_{B_1}|u_r|^p \eta^2\, dx \leq (p-1)\int_{B_1}|u_r|^p|\nabla
(r\eta)|^2\, dx
\end{equation}
holds for every radial $\eta\in C^1_c(B_1\setminus\{0\})$.
Since $\eta$ vanishes in a neighborhood of the origin and $u\in C^1$
away from the origin (see Remark~\ref{rem2:1}), we deduce by an 
approximation argument that \eqref{stab:cond} also holds for every
radial Lipschitz function $\eta$ vanishing on $\partial B_1$ and 
also in a neighborhood of the origin.

  We now prove that \eqref{stab:cond} also holds for every radial 
Lipschitz function vanishing on $\partial B_1$ ---but now 
not necessarily vanishing around $0$.
To see this, take $\zeta\in C^1(\mathbb{R}^n)$ such that $0\leq \zeta\leq 1$,
$\zeta\equiv 0$ in $B_1$, and $\zeta\equiv 1$ in $\mathbb{R}^n\setminus B_2$. 
Let $\zeta_{_\delta}(\cdot):=\zeta(\cdot/\delta)$ for every $\delta>0$.
Replacing $\eta$ by $\eta\zeta_\delta$ (which is Lipschitz and vanishes in a
neighborhood of $\{0\}$) in \eqref{stab:cond}, we obtain
\begin{equation}\label{lalala2}
(n-1)\int_{B_1}|u_r|^p\zeta_{_\delta}^2\eta^2\, dx\leq
(p-1)\int_{B_1}|u_r|^p|\nabla(r\eta\zeta_{_\delta})|^2\, dx.
\end{equation}
Now, using that $\eta$ and $\nabla\eta$ are bounded and
denoting $A_{\delta,2\delta}=B_{2\delta}\setminus \overline{B}_\delta$
we find
\begin{eqnarray}
 & &
 \hspace{-15mm}
\int_{B_1}|u_r|^p|\nabla(r\eta\zeta_{_\delta})|^2\, dx \nonumber \\
&=& \int_{B_1}|u_r|^p\Big\{|\nabla(r\eta)|^2\zeta_{_\delta}^2
+r^2\eta^2|\nabla\zeta_{_\delta}|^2
+\zeta_{_\delta}\nabla\zeta_{_\delta}\cdot\nabla(r^2\eta^2) \Big\}
dx\nonumber
\\
&\leq&\int_{B_1}|u_r|^p|\nabla(r\eta)|^2\zeta_{_\delta}^2\, dx
+C\int_{A_{\delta,2\delta}}|u_r|^p|\eta|\left\{\frac{r^2}{\delta^2}|\eta|
+\frac{r}{\delta}|\zeta_{_\delta}||\nabla(r\eta)|\right\}\, dx\nonumber \\
&\leq&\int_{B_1}|u_r|^p|\nabla(r\eta)|^2\zeta_{_\delta}^2\, dx
+C\int_{A_{\delta,2\delta}}|u_r|^p\, dx,\label{lalala3}
\end{eqnarray}
where $C$ denotes different positive constants. Since $u\in
W^{1,p}(B_1)$, the last term in \eqref{lalala3} tends
to zero as $\delta$ goes to zero. Therefore using that
$\zeta_\delta$ tend to 1 a.e. in $B_1$ as $\delta\to 0$,
\eqref{lalala3}, and \eqref{lalala2}, we obtain  by monotone
convergence that \eqref{stab:cond} holds for every radial
Lipschitz function $\eta$ vanishing on $\partial B_1$.


For $\alpha$ satisfying \eqref{eq:alpha} and $\varepsilon
\in(0,1)$, let
%
$$
\eta_\varepsilon (r) :=\left\{
\begin{array}{ll}
 \varepsilon ^{-\alpha }-1 & \text{ for }0\leq r\leq \varepsilon \\
r^{-\alpha }-1 & \text{ for }\varepsilon < r\leq 1,
\end{array}\right.
$$
%
a Lipschitz function vanishing on $\partial B_1$.
{From} inequality \eqref{stab:cond} applied with 
$\eta=\eta_\varepsilon$, we obtain
\begin{eqnarray*}
& &\hspace{-4mm} (n-1)\int_{B_{1}\setminus
B_\varepsilon}\hspace{-1mm}|u_r|^p (r^{-\alpha}-1)^2dx
+(n-1)(\varepsilon^{-\alpha}-1)^2
\int_{B_{\varepsilon}}|u_r|^pdx
\\
&
&\!\!\!\!\!\leq\hspace{-1mm}(p-1)\hspace{-1mm}\int_{B_{1}\setminus
B_\varepsilon}
\hspace{-2mm}|u_r|^p\left\{(1-\alpha)r^{-\alpha}\hspace{-1mm}-1\right\}^2dx
+(p-1)(\varepsilon^{-\alpha}-1)^2\int_{B_\varepsilon}|u_r|^pdx.
\end{eqnarray*}
Since $n\geq p$ it follows that
\[
(n-1)\int_{B_{1}\setminus B_\varepsilon}|u_r|^p
(r^{-\alpha}-1)^2dx \leq (p-1)\int_{B_{1}\setminus
B_\varepsilon}
\hspace{-0mm}|u_r|^p\left\{(1-\alpha)r^{-\alpha}-1\right\}^2dx.
\]
Developing the squares, using $n\geq p$ and \eqref{eq:alpha}, we
find the estimate
\begin{equation}\label{eq:6-19}
\int_{B_{1}\setminus B_\varepsilon}|u_r|^pr^{-2\alpha}dx\leq
\frac{C_{n,p}}{(n-1)-(\alpha-1)^2(p-1)} \int_{B_{1}\setminus
B_\varepsilon}|u_r|^pr^{-\alpha}dx.
\end{equation}
Throughout the proof\, $C_{n,p}$ (respectively $C_{n,p,\alpha}$)
denote different positive constants depending only on $n$ and $p$
(respectively on $n$, $p$, and $\alpha$). 

Now, choose a positive constant $C_{n,p,\alpha}$ such that
$$
\frac{C_{n,p}}{(n-1)-(\alpha-1)^2(p-1)}r^{-\alpha} \leq
\frac{1}{2}r^{-2\alpha}+C_{n,p,\alpha} \quad\text{for all
}r\in(0,1).
$$
Combining this inequality and \eqref{eq:6-19}, we obtain
\begin{equation*}
\int_{B_{1}\setminus B_\varepsilon}|u_r|^pr^{-2\alpha}dx \leq
C_{n,p,\alpha} \int_{B_{1}\setminus B_\varepsilon}|u_r|^pdx.
\end{equation*}
Now, let $\varepsilon\to 0$ to conclude
\begin{equation}\label{eq:6-20b}
\int_{B_{1}}|u_r|^pr^{-2\alpha}dx\leq C_{n,p,\alpha} 
\Vert  \nabla u \Vert_{L^{p}(B_1)}^p.
\end{equation}

In order to find a more precise expression on how
the previous constant $C_{n,p,\alpha}$ depends on $\alpha$, 
we apply \eqref{eq:6-20b} with the especial choice
$\alpha=\alpha_0$ given by
\[
\alpha_0=\frac{1}{2}+\frac{1}{2}\sqrt{\frac{n-1}{p-1}}
\in\left[1,1+\sqrt\frac{n-1}{p-1}\,\right).
\]
We deduce
\begin{equation}\label{eq:6-20c}
\int_{B_{1}}|u_r|^pr^{-1-\sqrt{\frac{n-1}{p-1}}}dx \leq
C_{n,p}\Vert \nabla u \Vert_{L^{p}(B_1)}^p.
\end{equation}
Finally, since $r^{-\alpha}\leq
r^{-\left(1+\sqrt{(n-1)/(p-1)}\right)}$ in $B_1$, \eqref{eq:6-19}
and \eqref{eq:6-20c} lead to the desired estimate \eqref{eq:6-11}
after letting $\varepsilon\to 0$.
\end{proof}

\section{Pointwise, $L^q$, and $W^{1,q}$  estimates}
\label{section3}

In this section we prove Theorem~\ref{thm1}.

\begin{proof}[Proof of Theorem {\rm \ref{thm1}}]
By \eqref{sobolev_rad}, $|u(1)|\le C_{n,p}\Vert u\Vert_\W1p$. Hence,
in view of the estimates that we need to prove, it suffices to
establish them with $u$ replaced by
$u-u(1)$ ---a positive solution vanishing on $\partial B_1$. 
Thus, through the proof, we assume
$$
u>0=u(1)\quad\text{in }B_1. 
$$

In case $n<p$, since $u\in \W1p$, the Sobolev embedding leads to
$u\in L^\infty(B_1)$ and $\Vert u\Vert_{L^\infty(B_1)}\le C_{n,p}
\Vert u\Vert_\W1p$.

In case $n\ge p$, let $\alpha$ satisfy \eqref{eq:alpha}. 
For $0<t<1$ we have
\begin{eqnarray}
u(t)
&=&\int _{t}^{1}-u_{r}r^{-(2\alpha-n+1)/p}
r^{(2\alpha-n+1)/p}dr\nonumber\\
&\le & C_{n,p}\left( \int
_{B_{1}}\hspace{-1mm}|u_{r}|^{p}
r^{-2\alpha}dx\right)^{\frac{1}{p}} \hspace{-1mm}\left(\int
_{t}^{1}\hspace{-1mm}
r^{p'(2\alpha-n+1)/p}dr\right)^{\frac{1}{p'}},
\label{eq:6-21}
\end{eqnarray}
by H\"older inequality and where $p'=p/(p-1)$. 
Using Lemma~\ref{lem2:plap} (which requires $n\ge p$) we deduce
\begin{equation}\label{eq:6-22}
u(t)\!\leq\!\frac{C_{n,p}\Vert u
\Vert_{W^{1,p}(B_1)} } {\left\{
(n-1)-(\alpha-1)^2(p-1)\right\}^{\frac{1}{p}}}\!
\left(\!\int_{t}^{1}
r^{p'(2\alpha-n+1)/p}dr\!\right)^{\frac{1}{p'}}\!\!\!
\end{equation}
for all $0<t<1$.

(a) Assume $n<p+4p/(p-1)$. By the remark made above, 
we may assume that $n\ge p$. Note that the integral in
\eqref{eq:6-22} is finite with $t=0$ whenever $p'(2\alpha-n+1)/p>-1$,
or equivalently 
\begin{equation}\label{eq:6-23}
\int_0^{1} r^{p'(2\alpha-n+1)/p}dr<+\infty\qquad
\textrm{if }\, \frac{n-p}{2}<\alpha.
\end{equation}
Recalling \eqref{rmk:cot_np} and since $n<p+4p/(p-1)$,
we can choose $\alpha$ (depending only on $n$ and $p$)
satisfying
$$
\frac{n-p}{2}<\alpha<1+\sqrt{\frac{n-1}{p-1}}.
$$
In addition, we may take $\alpha\ge 1$ ---as required in \eqref{eq:alpha}.
Now, the desired $L^\infty$ estimate follows from \eqref{eq:6-22} and
\eqref{eq:6-23}.

(b) Assume $n=p+4p/(p-1)$. Let $\varepsilon\in(0,1)$ and
\[
\alpha=1+\sqrt{\frac{n-1}{p-1}}-\varepsilon
=\frac{2p}{p-1}-\varepsilon.
\]
Since $\alpha$ satisfies \eqref{eq:alpha} in Lemma~\ref{lem2:plap},
\eqref{eq:6-22} yields
\begin{eqnarray}
u(t) &\leq& \frac{C_{p}\Vert u
\Vert_{W^{1,p}(B_1)}}{\varepsilon^{\frac{1}{p}}} \left(\int
_{t}^{1}r^{-1-2\varepsilon p'/p}dr\right)^{\frac{1}{p'}}
\nonumber\\
&\leq & 
\frac{C_{p}\Vert u
\Vert_{W^{1,p}(B_1)}}{\varepsilon}\ t^{-{2\varepsilon}/{p}},
\label{eq:6-24b}
\end{eqnarray}
for all $0<t < 1$ and $0<\varepsilon < 1$, where $C_p=C_{p,n}$
is a constant depending only  in $p$ (since here $n$ is a function
of $p$). It follows that $u\in L^q(B_1)$ for every
$1\leq q<\infty$.

In order to prove the pointwise estimate \eqref{eq:logpbd}, we
optimize the right-hand side of \eqref{eq:6-24b} with respect to
$\varepsilon$ by choosing $\varepsilon=\log{2}|\log t|^{-1}$.
Note that $\varepsilon$ is admissible since it belongs to $(0,1)$ 
if  $0<t<1/2$. With this choice of $\varepsilon$, \eqref{eq:6-24b}
yields
$$
u(t)\;\leq\; C_{p}\Vert u \Vert_{W^{1,p}(B_1)}|\log t| \qquad
\textrm{for}\quad 0<t<1/2.
$$
Using this and  that $u$ is positive and decreasing, 
the desired logarithmic estimate \eqref{eq:logpbd} follows.

(c) Assume $n>p+4p/(p-1)$ and $1\leq q<q_0$, for $q_0$ defined as in
\eqref{q0_and_q1}. For $\varepsilon\in (0,1)$, let
$$
\alpha=1+\sqrt{\frac{n-1}{p-1}}-\varepsilon.
$$
By   \eqref{eq:6-22}, we have
\begin{eqnarray}
u(t)& \leq & \frac{C_{n,p}\Vert u
\Vert_{W^{1,p}(B_1)}}{\varepsilon^{\frac{1}{p}}} \left(\int
_{t}^{1}
r^{\frac{p'}{p}\left(-n+2\sqrt{\textstyle\frac{n-1}{p-1}}
+3-2\varepsilon\right)}
dr\right)^{\frac{1}{p'}}\nonumber\\
& \leq & C_{n,p}\ \Vert u \Vert_{W^{1,p}(B_1)}
\frac{1}{\varepsilon^{\frac{1}{p}}}
t^{-\frac{1}{p}\left(n-2\sqrt{\textstyle\frac{n-1}{p-1}}-p-2
+2\varepsilon\right)},
\label{eq:6-25}
\end{eqnarray}
for $0<t<1$, where we have used that
$-\left(n-2\sqrt{{\frac{n-1}{p-1}}}-p-2\right)<0$ since
$n>p+4p/(p-1)$ ---see \eqref{rmk:cot_np}.

Now, from \eqref{eq:6-25} we obtain
$$
\int_{B_{1}}u^q\, dx \leq
 \frac{C_{n,p}^q
\Vert u \Vert_{W^{1,p}(B_1)}^q} 
{\varepsilon^\frac{q}{p}}
\int_0^1
t^{-\left(n-2\sqrt{\textstyle\frac{n-1}{p-1}}
-p-2+2\varepsilon\right)\frac{q}{p}}
t^{n-1}dt.
$$
If we set
$$
q=\frac{np}{n-2\sqrt{\frac{n-1}{p-1}}-p-2+3\varepsilon}<q_0
$$
and $\varepsilon>0$ is small enough, the second integral of the
previous inequality is finite. Hence, $u\in L^q(B_{1})$ for
every $1\leq q<q_0$. 

Finally, to prove the pointwise estimate \eqref{eq:estLqpws} we
consider \eqref{eq:6-25} and proceed as in part (b). Now, we need to
make $t^{-2\varepsilon/p}/\varepsilon^{1/p}$ small for given $t$.
We take $\varepsilon=\log 2|\log t|^{-1}$, which belongs to
$(0,1)$ if $0<t<1/2$. With this choice of $\varepsilon$,
\eqref{eq:6-25} leads to
$$
u(t)\;\leq\; C_{n,p}\ \Vert u \Vert_{W^{1,p}(B_1)} 
t^{-\frac{1}{p}\left(n-2\sqrt{\textstyle\frac{n-1}{p-1}}-p-2\right)}|\log
t|^\frac{1}{p}
$$
for $t\in (0,1/2)$. Since $u$ is positive and decreasing, this
leads to estimate \eqref{eq:estLqpws} in all $B_1$.

(d) Assume $g\ge 0$. We prove part (d1). 
Note that
\begin{equation}\label{16bis1}
\partial_r(r^{n-1}|u_r|^{p-1})=-\partial_r(r^{n-1}|u_r|^{p-2}u_r)
=r^{n-1}g(u)\ge 0\quad\text{in }B_1
\end{equation}
and hence $r^{n-1}|u_r|^{p-1}$ is a nonnegative and nondecreasing 
function of $r$. In particular
\begin{equation}\label{16bis2}
\Vert r^{n-1}|u_r|^{p-1}\Vert_{L^\infty(B_1)}\le |u_r(1)|^{p-1}.
\end{equation}
By \eqref{16bis1} and since $u$ is bounded in $B_1\setminus B_{1/2}$,
$r^{n-1}|u_r|^{p-1}$ is a $W^{1,1}(1/2,1)$ function of $r$ (indeed 
$W^{1,\infty}$). By the Sobolev embedding in one dimension, 
$r^{n-1}|u_r|^{p-1}$
is a continuous function up to $r=1$. Thus $|u_r(1)|$ in \eqref{16bis2}
is well defined and bounded. To control it, multiply equation 
\eqref{problem1} by $\varphi_\varepsilon(r):=\min\{1,\varepsilon^{-1}(1-r)\}$,
where $\varepsilon\in (0,1)$. This yields
\begin{equation}\label{16bis3}
\frac{1}{\varepsilon}\int_{B_1\setminus\overline{B}_{1-\varepsilon}}
|u_r|^{p-1}dx=\int_{B_1}g(u)\varphi_\varepsilon dx.
\end{equation}  
To fully justify \eqref{16bis3}, we need to multiply \eqref{problem1}
by $\varphi_\varepsilon\zeta_\delta$, where $\zeta_\delta$ vanishes 
around~$0$ as in the proof of Lemma~\ref{lem2:plap}. Then  we let $\delta\to 0$
and use that $u\in\W1p$ if $n\ge p$, and that $u\in C^1(\overline{B}_1)$ and 
$u_r(0)=0$ if $n<p$.

Now, letting $\varepsilon\to 0$ in \eqref{16bis3} we obtain 
$|u_r(1)|^{p-1}\le C_{n,p}\Vert g(u)\Vert_{L^1(B_1)}$. Thus, by 
\eqref{16bis2} we deduce that
\begin{equation}\label{16bis4}
\Vert r^{n-1}|u_r|^{p-1}\Vert_{L^\infty(B_1)}
\le C_{n,p} \Vert g(u)\Vert_{L^1(B_1)}.
\end{equation}

To control $\Vert \nabla u\Vert_{L^p(B_1)}$, assume first that $n<p$.
Then, 
\begin{eqnarray*}
\int_0^1r^{n-1}|u_r|^pdr 
& =&\int_0^1(|u_r|^{p-1}r^{n-1})^\frac{p}{p-1} r^{-\frac{n-1}{p-1}}dr\\
&\le&
\Vert r^{n-1}|u_r|^{p-1}\Vert^{\frac{p}{p-1}}_{L^\infty(B_1)}
\int_0^1r^{-\frac{n-1}{p-1}}dr.
\end{eqnarray*}
Since the last integral is finite, this and \eqref{16bis4} lead to 
\eqref{eq:estGLq}. 

In case $n\ge p$, we use Lemma~\ref{lem2:plap}.
We take $\alpha$  satisfying \eqref{eq:alpha} and depending 
only on $n$ and $p$. The lemma gives that 
\begin{eqnarray}\label{esss1}
&& \hspace{-12mm}\int_{B_1}r^{-2\alpha}|u_r|^pdx \\ \label{esss2}
&\le &C_{n,p}
\int_{B_1}|u_r|^pdx  = C_{n,p} \int_{B_{r_0}}|u_r|^pdx 
+C_{n,p} \int_{B_1\setminus \overline{B}_{r_0}}|u_r|^pdx ;
\end{eqnarray}
here we choose $r_0\in (0,1)$ satisfying $r_0^{-2\alpha}\geq 2C_{n,p}$.
Since $C_{n,p}\int_{B_{r_0}}|u_r|^pdx
\leq (1/2) \int_{B_{r_0}}r^{-2\alpha}|u_r|^pdx
\leq (1/2) \int_{B_1}r^{-2\alpha}|u_r|^pdx$, we can absorb the first
term in the right-hand side of \eqref{esss2} into \eqref{esss1},
and deduce that
\begin{equation}\label{16bis5}
\int_{B_{1}}|u_r|^pdx\le\int_{B_{1}}r^{-2\alpha}|u_r|^pdx\le 2C_{n,p}
\int_{B_1\setminus\overline{B}_{r_0}}|u_r|^pdx.
\end{equation}
Note that $r_0$ depends only on $n$ and $p$, 
and thus, since $u$ is decreasing, $u(r_0)^{p-1}\le C_{n,p}\Vert u^{p-1}
\Vert_{L^1(B_{r_0})}$. Using this and \eqref{16bis4}, we have
\begin{eqnarray*}
&&\hspace{-12mm}\int_{B_1\setminus\overline{B}_{r_0}}|u_r|^pdx
= C_n \int_{r_0}^1|u_r|^pr^{n-1}dr\\
&\le & C_n \Vert r^{n-1}|u_r|^{p-1}\Vert_{L^\infty(B_1)}\int_{r_0}^1-u_rdr
\le  C_{n,p}\Vert g(u)\Vert_{L^1(B_1)}
\Vert u^{p-1}\Vert_{L^1(B_1)}^\frac{1}{p-1}.
\end{eqnarray*} 
After using Young's inequality, this bound in 
$B_1\setminus\overline{B}_{r_0}$ and estimate 
\eqref{16bis5} yield the desired bound 
\eqref{eq:estGLq} in all $B_1$.

Next, we prove (d2). We consider first the case $n<p+4p/(p-1)$. 
Observe that by part (a) we have that $u\in L^\infty(B_1)$.
Hence, applying the regularity theory for the $p$-Laplacian
(see \cite{Lie88}) it follows that $u\in C^{1,\beta}(\overline{B}_1)$ for some
$0<\beta<1$. Moreover, by the $L^\infty$ estimate of
part~(a) above and estimate \eqref{eq:estGLq}, there exists a constant 
depending on $n$, $p$, $\Vert u\Vert_{L^1(B_1)}$, and upper bounds in $g$
such that  estimate \eqref{eq:Wbd} holds. This concludes the proof
for $n<p+4p/(p-1)$.

Assume $n\geq p+4p/(p-1)$. We establish parts (d2) and (d3) at the same time. 
Observe that it is enough to prove our estimates in $B_{1/4}$. Indeed, 
since $u$ is decreasing, we have bounds for the supremum of $u$ in 
$B_1\setminus B_{1/5}$ in terms of $\Vert u\Vert_{L^1(B_1)}$.
Thus by \cite{Lie88} we have that 
$u\in W^{1,q}(B_1\setminus \overline{B}_{1/5})$
for all $q<\infty$. Moreover, $\Vert u\Vert_{W^{1,q}(B_{1}\setminus
\overline{B}_{1/5})}\le C$, where $C$ is a constant depending on $n$, $p$,
$\Vert u\Vert_{L^1(B_1)}$, and upper bounds in $g$.

We choose $\tilde{\rho}\in (1/4,1/2)$ such that
\begin{equation}\label{eq:bdrho}
0<-u_r(\tilde\rho)=\frac{u(1/4)-u(1/2)}{1/4}\le
4\|u\|_{L^\infty(1/4,1)}\le C_{n,p}\Vert u \Vert_{W^{1,p}(B_1)},
\end{equation}
where in the last inequality we have used \eqref{sobolev_rad}.
On the other hand, from \eqref{problem1:radial} we have
\begin{equation}\label{eq:6-27}
\partial_r(|u_r|^{p-2}u_r)=-\frac{n-1}{r}|u_r|^{p-2}u_r-g(u)
\le -\frac{n-1}{r}|u_r|^{p-2}u_r,
\end{equation}
for $r\in(0,1)$. We integrate \eqref{eq:6-27} with respect to $r$,
from $t$ to $\tilde\rho$, and use \eqref{eq:bdrho}, to obtain
\begin{eqnarray}
-|u_r(t)|^{p-2}u_r(t)&\leq& -|u_r(\tilde\rho)|^{p-2}u_r(\tilde\rho)
-(n-1)\int_t^{\tilde\rho}\frac{|u_r|}{r}^{p-2}u_r\ dr
\nonumber \\
&\leq & C_{n,p}\Vert u \Vert_{W^{1,p}(B_1)}^{p-1}+
(n-1)\int_t^{1/2}\frac{|u_r|}{r}^{p-1}dr, \label{eq:6-28}
\end{eqnarray}
for all $0<t<1/4$. Next, we estimate the last integral in 
\eqref{eq:6-28} using H\"older inequality, to find
\begin{eqnarray}
\int_t^{1/2}\frac{|u_r|}{r}^{p-1}dr&
=&\int_t^{1/2}|u_r|^{p-1}r^{-(2\alpha-n+1)/p'}
r^{(2\alpha-n+1-p')/p'}dr \nonumber
\\
&\hspace{-4cm}\leq &\hspace{-2.2cm}
C_n\left(\int_{B_{1/2}}|u_r|^pr^{-2\alpha}dx\right)^\frac{1}{p'}
\left(\int_t^{1/2} r^{p(2\alpha-n+1-p')/p'}dr\right)^\frac{1}{p}.
 \label{eq:holdplap}
\end{eqnarray}
Let $\varepsilon\in (0,1)$ and
$$
\alpha=1+\sqrt\frac{n-1}{p-1}-\varepsilon.
$$
Applying Lemma~\ref{lem2:plap} in \eqref{eq:holdplap} and
using that $p(2\alpha-n+1-p')+p'<0$ since $n\geq p+4p/(p-1)$, we
deduce from \eqref{eq:6-28} and \eqref{eq:holdplap} the following
estimate
\begin{equation}\label{eq:6-29}
-u_r(t)\leq \frac{C_{n,p}}{\ve^{1/p}}\ \Vert u \Vert_{W^{1,p}(B_1)}
\
t^{-\frac{1}{p}\left(n-2\sqrt{\textstyle\frac{n-1}{p-1}}-2+2\varepsilon\right)},
\end{equation}
for $0<t<1/4$.

Now we use \eqref{eq:6-29} to obtain, for $q\ge 1$,
\begin{equation}\label{eq:6-30}
\int_{B_{1/4}}\hspace{-1mm}|u_r|^qdx\leq \frac{C_{n,p}}{\ve^{q/p}}
\Big(\Vert u \Vert_{W^{1,p}(B_1)} \Big)^q\hspace{-0.2mm}
\int_0^1\hspace{-0.8mm}
r^{n-1-\frac{q}{p}\left(n-2\sqrt{\textstyle\frac{n-1}{p-1}}-2
+2\varepsilon \right)}dr.
\end{equation}
If we set
$$
q=\frac{np}{n-2\sqrt{\frac{n-1}{p-1}}-2+3\varepsilon}<q_1
$$
the second integral in \eqref{eq:6-30} is finite for every
$\varepsilon>0$. Hence, $u_r\in L^q(B_{1/4})$ for every $q<q_1$. The
estimate \eqref{eq:Wbd} follows, using in addition estimate \eqref{eq:estGLq}.

Finally, the pointwise estimate \eqref{eq:pdWdb}
follows from \eqref{eq:6-29}
by choosing $\ve=\log 4|\log t|^{-1}$ for $0<t<1/4$.
\end{proof}


\section{The extremal solution}
\label{section4}

This section is devoted to prove Theorem~\ref{thm2}. 
By estimate \eqref{eq:estGLq}
of Theorem~\ref{thm1}, now our task
is to bound $u_\lambda^{p-1}$ and $f(u_\lambda)$ in $L^1(B_1)$ uniformly
in $\lambda\in (0,\lambda^*)$. We prove this fact using the growth condition 
\eqref{superlinear} on $f$ together with the radially decreasing property 
of the minimal solutions $u_\lambda$. These bounds lead to a control of
$u_\lambda$ in $\W1p$ uniformly in $\lambda$, by estimate \eqref{eq:estGLq}.

\begin{proof}[Proof of Theorem {\rm \ref{thm2}}]
We start noting that, for $\lambda\in (0,\lambda^*)$, the minimal
solution $u_\lambda\in L^\infty(B_1)$ is radially decreasing. 
This follows from a general result of \cite{DS04} on radially decreasing 
symmetry, which only requires $f$ to be positive in $(0,\infty)$ 
and locally Lipschitz in $[0,\infty)$ ---see Corollary~1.1 in \cite{DS04}. 
We can give, however, a simple proof in our situation of minimal solutions, 
as follows. First, the minimal solution $u_\lambda$ can be constructed by 
monotone iteration (see \cite{CS05}), starting with $u^0\equiv 0$, by
solving $-\Delta_pu^k=\lambda f(u^{k-1}(r))$. Since the right-hand side is given 
and radial, the unique solution $u^k$ must be radial. Thus, the limit 
$u_\lambda$,
as $k\to \infty$, is radial. Now, integrating our equation 
$-\partial_r(r^{n-1}|\partial_ru_\lambda|^{p-2}\partial_ru_\lambda)
=r^{n-1}\lambda f(u_\lambda)>0$ in $r$, from $0$ to $t\in (0,1)$, we
get that $\partial_ru_\lambda(t)<0$. 

For $\lambda\in (0,\lambda^*)$, let $\rho_\lambda\in (1/2,1)$ be such that
\begin{eqnarray*}
|\partial_ru_\lambda(\rho_\lambda)|^{p-1}
&=&(-\partial_ru_\lambda(\rho_\lambda))^{p-1}
=\left(\frac{u_\lambda(1/2)-u_\lambda(1)}{1/2}\right)^{p-1}\\
&=& (2u_\lambda(1/2))^{p-1}\le C_{n,p}\Vert u_\lambda^{p-1}\Vert_{L^1(B_{1/2})},
\end{eqnarray*}
where we have used that $u_\lambda$ is radially decreasing. Since 
$r^{n-1}|\partial_ru_\lambda|^{p-1}$ is increasing in $r\in (0,1)$ 
by \eqref{16bis1},
the previous estimate yields
\begin{equation}\label{19bis1}
\Vert r^{n-1}|\partial_ru_\lambda|^{p-1}\Vert_{L^\infty(B_{1/2})}
\le \rho_\lambda^{n-1}|\partial_ru_\lambda(\rho_\lambda)|^{p-1}
\le C_{n,p}\Vert u_\lambda^{p-1}\Vert_{L^1(B_{1/2})}.
\end{equation}
Multiply equation (\extp) by $\psi(r)=\min\{1,4(1/2-r)^+\}$,
to obtain
\begin{eqnarray}
\Vert \lambda f(u_\lambda)\Vert_{L^1(B_{1/4})}
&\le&C_{n,p}
\int_{1/4}^{1/2}r^{n-1}|\partial_ru_\lambda|^{p-1}dr
\nonumber\\
&\le & C_{n,p}\Vert u_\lambda^{p-1}\Vert_{L^1(B_{1/2})},
\label{19bis2}
\end{eqnarray}
where we have used \eqref{19bis1}. Note that here we can test
the equation with $\psi$, a function which does not 
vanish around the origin, since we know that $u_\lambda$ is an
energy solution in all of $B_1$.

Next, we use assumption \eqref{superlinear} on $f$ to ensure, given 
any $\delta>0$, that $\lambda f(t)\ge \frac{1}{\delta}t^{p-1}-C_\delta$
for all $t>0$ and $\lambda\in (\lambda^*/2,\lambda^*)$, where $C_\delta$ 
does not depend on $\lambda$. This combined with \eqref{19bis2} leads to
\begin{equation}\label{19bis3}
\Vert u_\lambda^{p-1}\Vert_{L^1(B_{1/4})}
\le C_{n,p}\delta\Vert u_\lambda^{p-1}\Vert_{L^1(B_{1/2})}+C_\delta.
\end{equation}
Now, since $u_\lambda$ is decreasing in $r$, we have
\begin{eqnarray}\label{19bis4}
\Vert u_\lambda^{p-1}\Vert_{L^1(B_{1/2}\setminus\overline{B}_{1/4})}
\le C_{n,p} u_\lambda^{p-1}(1/4)
\le C_{n,p} \Vert u_\lambda^{p-1}\Vert_{L^1(B_{1/4})}.
\end{eqnarray}
 This combined with \eqref{19bis3} yields, after taking $\delta$ small 
enough, $\Vert u_\lambda^{p-1}\Vert_{L^1(B_{1/4})}\le C$ for a constant
$C$ independent of $\lambda$. This bound and the argument used in 
\eqref{19bis4} done now on $B_1\setminus\overline{B}_{1/4}$ leads to a 
uniform in $\lambda$ bound for $\Vert u_\lambda^{p-1}\Vert_{L^1(B_1)}$.

The previous bound gives also a control for  
$\Vert f(u_\lambda)\Vert_{L^1(B_{1/4})}$,
by \eqref{19bis2}. Since $f$ is increasing, $f(u_\lambda)$ is decreasing
in $r$ and the argument used above allows to control 
$\Vert f(u_\lambda)\Vert_{L^1(B_1)}$ uniformly in $\lambda$. Thus,
$$
\Vert u_\lambda^{p-1}\Vert_{L^1(B_1)} +\Vert f(u_\lambda)\Vert_{L^1(B_1)}
\le C
$$  
for some constant $C$ independent of $\lambda\in (0,\lambda^*)$. By estimate 
\eqref{eq:estGLq} of Theorem~\ref{thm1}, we deduce a 
bound for $\Vert \nabla u_\lambda\Vert_{L^p(B_1)}$ independent of $\lambda$. 
Since $u_\lambda\big|_{\partial B_1}\equiv 0$, this yields a bound for 
$\Vert u_\lambda\Vert_\W1p$.

Thus, we have that there exists $v\in
W^{1,p}_0(B_1)$ and some subsequence still denoted by $u_\lambda$
such that, as $\lambda\to \lambda^*$, $u_\lambda\rightharpoonup
v$  weakly in $W^{1,p}_0(B_1)$, $u_\lambda \rightarrow v$ strongly in
$L^p(B_1)$, and $u_\lambda\rightarrow v$ a.e. in $B_1$. 
Since $u_\lambda$ is increasing in $\lambda$, every sequence $u_\lambda(r)$
tends as $\lambda\to\lambda^*$ to $u^*(r)$,
by the definition  \eqref{extremal} of $u^*$. Hence,
$v=u^*\in W^{1,p}_0(B_1)$.

Next, we show that $u^*$ is an
energy solution of (\exst). We want to
pass to the limit in
\begin{equation}\label{distributional2}
\int_{B_1}|\nabla u_\lambda|^{p-2}\nabla
u_\lambda\cdot\nabla\varphi\, dx=\int_{B_1}\lambda
f(u_\lambda)\varphi\, dx\quad \textrm{for }\varphi\in C_c^1(B_1).
\end{equation}
Since $f(u_\lambda)$ increases to $f(u^*)$ and 
$\Vert f(u_\lambda)\Vert_{L^1(B_1)}$ is uniformly bounded in $\lambda$,
the monotone convergence theorem gives that $f(u^*)\in L^1(B_1)$ and that
the limit as $\lambda\to\lambda^*$ of the right-hand side in 
\eqref{distributional2} is
$$
\lambda^* \int_{B_1} f(u^*)\varphi\, dx.
$$
To pass to the limit in the left-hand side of \eqref{distributional2}, note 
first that we may assume that $\varphi$ is radial ---simply by integrating 
first with respect to the ``angles'' $\theta\in \partial B_1$ in 
\eqref{distributional2}. Now, the left-hand side of \eqref{distributional2},
up to a multiplicative constant, reads
\begin{equation}\label{19bis5}
\int_0^1r^{n-1}|\partial_ru_\lambda|^{p-1}(-\partial_r\varphi)(r)dr.
\end{equation}  
Using \eqref{16bis1} and that the $L^1(0,1)$ norm of $r^{n-1}f(u_\lambda(r))$
is bounded uniformly in~$\lambda$, we deduce that 
$r^{n-1}|\partial_ru_\lambda|^{p-1}$ is
bounded in $W^{1,1}(0,1)$ ---the Sobolev space in one dimension---
uniformly in $\lambda\in (0,\lambda^*)$. Since this space is compactly
embedded in $L^1(0,1)$, we conclude that $r^{n-1}|\partial_r u_\lambda|^{p-1}$ 
converges strongly in $L^1(0,1)$ to $r^{n-1}|\partial_ru^*|^{p-1}$.
Hence, we can pass to the limit in \eqref{19bis5} and conclude   
$$
\int_{B_1} |\nabla u^*|^{p-2}\nabla u^*\cdot\nabla \varphi\,
dx\displaystyle = \lambda^* \int_{B_1} f(u^*)\varphi\, dx,
$$
for every $\varphi\in C_c^1(B_1)$, proving the claim.

Clearly, $u^*$ is a radially symmetric solution and it is nonincreasing. 
Hence, by regularity, $u^*\in
C^1(\overline{B}_1\setminus\{0\})$ since $u^*\in
L^\infty(B_1\setminus\{0\})$. Therefore, using that 
$u^*\in W^{1,p}_0(B_1)$ and $f>0$, the argument
at the end of the first paragraph of this proof ---or Hopf Lemma 
(see \cite{Peral} or Lemma A.3 in \cite{Sk})---, 
lead to $u_r^*(t)<0$ for all $t\in(0,1)$. 
That is, $u^*$ is radially decreasing.

At this point we are able to prove that $u^*$ is semi-stable in the
sense of Definition~\ref{def:stab}. That is,
\begin{equation}
\int_{B_1}\lambda^*f'(u^*)\xi^2dx\le 
\int_{B_1}(p-1)|u_r^*|^{p-2}|\xi_r|^2 dx
\end{equation}
for every radially symmetric function $\xi\in C^1_c(\b1z)$. This
follows by passing to the limit as $\lambda\to\lambda^*$ in
the corresponding semi-stability property for the minimal solutions
$u_\lambda$. On the left-hand sides, we simply, use Fatou's lemma 
(recall that $f'\ge 0$). On the right-hand side, we use that $\xi$ has
 compact support in $\b1z$, and that in such compact set we have $C^{2,\beta}$
estimates for $u_\lambda$ uniformly in $\lambda\in (0,\lambda^*)$. 
This holds since $|\partial_ru_\lambda^*|>0$ in such compact set,
and the same is true for $|\partial_ru_\lambda|$ uniformly in 
$\lambda\in (0,\lambda^*)$, and thus we deal with uniformly 
elliptic equations. Thus $|\partial_ru_\lambda|$ converges uniformly
to $|u^*_r|$, and the same is true for the quantities 
$|\partial_r u_\lambda|^{p-2}$, since $|u^*_r|> 0$ ---note that 
the exponent $p-2$ could be negative.   

Finally, the regularity statements in the theorem
follow as a consequence of Theorem \ref{thm1}.
\end{proof}

\centerline{{\bf Acknowledgments}}

\smallskip
The authors were supported by the Ministerio de Educaci\'on y
Ciencia (Spain), grant MTM2005-07660-C02-01. The  authors are
partially sponsored by the European Science Foundation (ESF) PESC
Programme ``Global''.

The second author is also supported by the E.U. RTN program
``MULTIMAT''  MRTN-CT-2004-505226.

The third author was also supported by CMUC/FCT (Coimbra, Portugal).
\medskip


\end{document}